\documentclass[10pt]{article}     
\usepackage{amssymb} 
\usepackage{amsmath,amscd} 
\newtheorem{thm}{Theorem}[section]    
\newtheorem{lem}[thm]{Lemma}          
\newtheorem{proposition}[thm]{Proposition}

\newtheorem{defn}{Definition}

\title{Rational formality of function spaces  } 
\date{}                   
\author{Micheline VIGU\'E-POIRRIER                    
    }                
\begin{document}

\maketitle 

\begin{abstract}   
Let $X$ be a nilpotent space such that there exists $N\geq 1$ with $H^N(X,\mathbb Q) \ne 0$ 
and $H^n(X,\mathbb Q)=0$ if $n>N$. Let $Y$ be a m-connected space with $m\geq N+1$ and $H^*(Y,\mathbb Q)$ is finitely generated as algebra. We assume that  the odd part of the rational Hurewicz homomorphism: $\pi _{odd}(X)\otimes \mathbb Q\rightarrow H_{odd}(X,\mathbb Q)$ is non-zero. We prove that if the space $\mathcal F(X,Y)$ of continuous maps from $X$ to $Y$ is rationally formal, then $Y$ has the rational homotopy type of a finite product  of Eilenberg Mac Lane spaces. At the opposite, we 
exhibit an example of a  rationally formal space $\mathcal F(S^2,Y)$ where $Y$ is not rationally equivalent to a product  of Eilenberg Mac Lane spaces.

\end{abstract}
\vspace{5mm}\noindent {\bf AMS classification}: 55P62, 55P35

\vspace{2mm}\noindent {\bf Key words}: function spaces, formal space, Sullivan model, Lie model.

\section{Introduction}
All the spaces we consider, have the rational homotopy type of  CW complexes of finite type. Let $X$ be a nilpotent space
 such that there exists $N\geq 1$ with $H^N(X,\mathbb Q) \ne 0$ and $H^n(X,\mathbb Q)=0$ if $n>N$. Let $Y$ a m-connected 
space with $m\geq N+1$.\\
Under these hypothesis, the space $\mathcal F(X,Y)$ of continuous maps from $X$ to $Y$ is simply connected. The rational homotopy type of this space has been determined in \cite{BS},  \cite{Hae}, \cite{Si}, \cite{Vi} where Sullivan models or Lie models are computed.

In his paper {\it L'homologie des espaces fonctionnels} \cite{Th}, Thom studied the homotopy type of the space of 
continous maps from $X$ to $Y$ homotopic to a given map $f$. He proved that if $Y$ is an Eilenberg-Mac Lane space 
then $\mathcal F(X, Y)$ has the homotopy type of a product of Eilenberg Mac Lane spaces. This implies that 
if $H^*(Y,\mathbb Q)$ is a free commutative algebra, then $H^*(\mathcal F(X,Y),\mathbb Q)$ is a free commutative algebra 
for any  $X$. Another proof is given in \cite{Vi}. Recall that a 1- connected space has the rational
 homotopy of a product of Eilenberg-Mac Lane spaces if and only if its cohomology algebra is free commutative.

In rational homotopy theory, the notion of rational formality plays a crucial role (see below 2.4), since the rational homotopy 
type of a formal space is entirely determined by the data of the singular cohomology algebra. \\
In the following, we will write "formality" instead of "rational formality" since we always work with nilpotent spaces in the context of rational homotopy theory.

\medskip

An open question
 is a converse to Thom's result.\\
{\bf Question:} What conditions  should be satisfied by $X$ and $Y$ if $\mathcal F(X,Y)$ is a formal space?
\medskip

In \cite{DV} it is proved that for $X=S^1$ and $H^*(Y,\mathbb Q)$  a finitely generated algebra,  
if $\mathcal F(X,Y)$ is a formal space, then $Y$ has the rational homotopy type of a product of Eilenberg Mac Lane spaces. The proof relies on the theory 
of Sullivan minimal models. With similar methods, Yamaguchi proves in \cite{Y} that if $Y$ satisfies
 $\dim H^*(Y,\mathbb Q)<+\infty $ and $\dim \pi _*(Y)\otimes \mathbb Q<+\infty $, then the formality of $\mathcal F(X,Y)$
 implies that  $Y$ has the rational homotopy type of a product of odd dimensional spheres. In this paper we prove:

\bigskip

{\bf Main Theorem :} {\sl Let $X$ be a nilpotent space such that there exists $N\geq 1$ with $H^N(X,\mathbb Q) \ne 0$ 
and $H^n(X,\mathbb Q)=0$ if $n>N$. Let $Y$ be a m-connected space with $m\geq N+1$. We assume that:  
\begin{enumerate}
\item the odd part of the rational Hurewicz homomorphism: $\pi _{odd}(X)\otimes \mathbb Q\rightarrow H_{odd}(X,\mathbb Q)$ is non-zero.
\item $H^*(Y,\mathbb Q)$ is finitely generated as algebra.
\end{enumerate}
Then, if $\mathcal F(X,Y)$ is formal, $Y$ has the rational homotopy type of a finite product of Eilenberg Mac Lane spaces.
}
\medskip

{\bf Corollary}: Under the hypothesis of the Main Theorem, $\mathcal F(X,Y)$ has the rational homotopy type of a product of Eilenberg Mac Lane spaces.

\medskip

{\bf Remark 1:} If $\pi_1 (X)$ is non-zero, assumption 1. is satisfied. If $\pi _1(X)=0$, the dual of the rational Hurewicz map can be identified with the map: $H^*(\bigwedge V,d)\rightarrow V$ induced by the projection where $(\bigwedge V,d)$ is a minimal Sullivan model of $X$, (\cite{FHT}, page 210). So assumption 1. is equivalent to the following: {\sl $X$ has a minimal Sullivan model $(\bigwedge V,d)$ with ${\rm Ker} d\cap V_{odd}\ne0$}. 

{\bf Remark 2:} Suppose $X$ is formal and there exists $q$ odd such that \newline $H^q(X,\mathbb Q) \ne 0$. Let $2d+1= \inf \lbrace q, \, \, {\rm odd}, \, \,  H^q(X,\mathbb Q) \ne 0\rbrace$, then there exists a nonzero element $a\in H^{2d+1}(X)$ and $a$ does not belong to $H^+(X)\cdot H^+(X)$. So $X$ has a minimal bigraded model in the sense of \cite{HS}, $  \rho:  (\bigwedge V,d)\rightarrow H^*(X)$ with a generator $t\in V_0$, $(dt=0)$, $\vert t\vert=2d+1$  and $\rho (t)=a$. Such a space satisfies assumption 1. of the theorem.\\
Example 6.5 in \cite{HS} provides a non formal space $X$ satisfying the  assumption 1. of the main theorem.

{\bf Remark 3:} Probably assumption 2. is not necessary.
\bigskip

The proof uses simultaneously the theory of minimal Quillen models of a space in the category of Lie differential graded
 algebras and the theory of minimal Sullivan models of a space in the category of commutative differential graded
 algebras. For that reason, we should ask the connectivity hypothesis on $Y$ to ensure $\mathcal F(X,Y)$ to be $1$-connected.
The idea of the proof is the following: we use Lie models to prove that, under the hypothesis of the theorem, the formality of $\mathcal F(X,Y)$ implies the formality of some $\mathcal F(S^p,Y)$ with $p$ odd (theorem 3.1). Then we work with a Sullivan model for $\mathcal F(S^p,Y)$, $p$ odd,   and we mimick the proof of \cite{DV} (theorem 3.5).

In the last section, we give an explicit example where $X=S^2$, $Y$ is not a product of Eilenberg Mac Lane spaces and however $\mathcal F(S^2,Y)$ is formal. This proves that assumption 1. is necessary.

\section{Algebraic models  in rational homotopy theory}

All the graded vector spaces, algebras, coalgebras and Lie
algebras $V$ are defined over $\mathbb Q$ and are supposed of
finite type, i.e. dim $V_n <\infty$ for all $n$.\\
If $v$ has degree $n$, we denote $\vert v\vert =n$.

For precise definitions we refer to \cite{FHT} or \cite{Hal}.\\
If $V =\{V_i\}_{i \in \mathbb Z}  $ is a (lower)
 graded $\mathbb Q$-vector space  (when we need upper graded vector space
we put
 $V_i= V^{-i}$ as usual.)

We denote by $sV$  the suspension of $V$, we have:
$(sV)_n =V_{n-1}\,, \,\,
(sV)^n =V^{n+1}. $

A morphism betwen two differential graded vector spaces is called a quasi-isomorphism 
if it induces an isomorphism in homology.

\subsection{Commutative differential graded algebras}
In the following we consider only commutative differential algebras graded in positive degrees, $(A=\oplus_ {n\geq 0} A^n,d)$ with a differential $d$
 of degree $+1$  satisfying $H^0(A,d)=\mathbb Q$ and $\dim A^n$ is finite for all $n$.  We denote by CDGA the category of commutative differential graded algebras. Such an algebra is called a commutative cochain algebra.
 If $V =\{V^i\}_{i\geq 1} $ is a
 graded $\mathbb Q$-vector space    we denote by $ \bigwedge V$
 the free graded commutative algebra generated by $V$. A commutative cochain algebra of the form $(\bigwedge V,d)$ where $d$ satifies some nilpotent conditions  is called a Sullivan algebra,(\cite{FHT},12) . A Sullivan algebra is called minimal if $dV\subset  \bigwedge^+V.\bigwedge^+V$.
 \begin{defn} 
A Sullivan model for a commutative cochain algebra $(A,d)$ is a quasi-isomorphism of differential graded algebras:
$$(\bigwedge V,d)\rightarrow (A,d)$$
with $(\bigwedge V,d)$ a Sullivan algebra.

If $d$ is minimal, we say that $(\bigwedge V,d)$ is a minimal Sullivan model.
\end{defn}
Any commutative cochain algebra has a minimal Sullivan model. If $H^1(A,d)=0$, then two minimal Sullivan models are isomorphic.
Any path connected
space $X$ admits a Sullivan  model which is the Sullivan model of the cochain algebra $A_{PL}(X)$,  where $A_{PL}$ denotes the contravariant functor of piecewise linear differential forms.
 Any simply connected space admits a minimal Sullivan model
 which contains all the informations on the rational
homotopy type of the space, (\cite{FHT}, §12).

\medskip
\begin {proposition}
Let $(\bigwedge V,d)$ be a Sullivan algebra such that $\dim V^i<\infty $ for all $i$ and $\dim H^*((\bigwedge V,d)<\infty $. Then there exist a commutative cochain algebra $(A,d)$ with $\dim A<\infty $ and a quasi-isomorphism of differential graded algebras: $(\bigwedge V,d)\rightarrow (A,d)$. 
\end{proposition}

{\bf Proof.} Let $p$ be an integer such that $H^n(\bigwedge V,d)=0$ for all $n>p$. We define a graded subspace $I\subset (\bigwedge V)$ so that:
$$I^k=0, \, \, k<p, \quad I^k=(\bigwedge V)^k, \, \, k>p$$
$$I^p\oplus (ker d)^p=(\bigwedge V)^p$$
Then $I=\oplus _k I^k$ is a differential ideal and $H^k(I,d)=0$ for all $k$. Put $A=(\bigwedge V/I, d)$, then $\dim A<\infty $  and the projection: $(\bigwedge V,d)\rightarrow (A,d)$  is a quasi-isomorphism. 
\medskip 

\subsection{Differential graded Lie algebras}
In the following we consider only differential graded Lie algebras: $L=(L_i)_{i\geq 1}$ and the differential  has degree $-1$.

Recall that $TV$ denotes the tensor algebra on a graded vector space $V$, it is a graded Lie algebra if we endow it with
 the commutator bracket. The sub Lie algebra generated by $V$ is called the free graded Lie algebra
 on $V$ and it is denoted $\mathbb L(V)$. A free differential Lie algebra $(\mathbb L(V),\partial)$ is called minimal 
if $\partial(V)\subset [\mathbb L(V),\mathbb L(V)]$.\\
\begin{defn}
 A free Lie model of a chain Lie algebra $(L,d)$ is a quasi-isomorphism of differential Lie algebras of the form
$$m: (\mathbb L(V),\partial)\rightarrow (L,d)$$
\end{defn}
If $\partial $ is minimal, it is called a minimal free Lie model.
Every chain Lie algebra $(L,d)$ admits a minimal free Lie model, unique up to isomorphism.
Every simply connected space has a minimal free Lie model 
$(\mathbb L(V),\partial)$ containing all the informations on the rational homotopy type of the space, (\cite{FHT}, §24)
 called the minimal Quillen model. 

A differential graded Lie algebra is called a model for a space $Y$ if its minimal free Lie model is the minimal Quillen of the space.

\subsection{ Dictionary between Sullivan models and Lie models}
A way of constructing Sullivan algebras from  differential graded Lie algebras is given by
the  functor $\mathcal C^*$ which is obtained by dualizing the Cartan-Chevalley construction that associates a cocommutative differential coalgebra to a differential Lie algebra,(\cite{FHT}, §23 ). In fact $C^*(L,d_L)$ is a Sullivan algebra $(\bigwedge V,d)$ with differential $d=d_0+d_1$,  $d_0(V)\subset V$ and $d_1(V)\subset \bigwedge^2(V)$.
More precisely $V$ and $sL$ are dual graded vector spaces, $d_0$ is dual of $d_L$, $d_1$ corresponds by duality to the Lie bracket on $L$.

\subsection{Formal commutative differential algebras and formal spaces}
\begin{defn}
A commutative cochain algebra $(A,d)$ is {\bf formal}  if its minimal model is quasi-isomorphic to $(H=H^*(A,d),0)$.
A space $M$ whose Sullivan minimal model
$(\bigwedge V,d)$ is quasi-isomorphic to $(H^\ast
(M),0)$ is called formal.
\end{defn}
 Examples  of formal spaces are given by  Eilenberg-Mac Lane spaces,
spheres, complex projective spaces.
Connected compact K\"ahler manifolds (\cite{DGMS}) and  quotients of compact
connected Lie groups by  closed subgroups  of the same rank are formal. Symplectic manifolds need not be formal.
Product and wedge of formal spaces are formal.

We  give now a  property for the conservation of formality between two cochain algebras, it will be a key point in the proof of the main theorem. A variant of this result is proved in \cite {DV}.

\begin{proposition}
Let $(A,d_A)$ and $(B,d_B)$ two commutative differential graded algebras satisfying $H^1(A)=0$. We assume that there exist two CDGA morphisms $f: (A,d_A)\rightarrow (B,d_B)$ and $g: (B,d_B)\rightarrow (A,d_A)$ satisfying: $g\circ f=Id$. If $(B,d_B)$ is formal, then $(A,d_A)$ is formal.
\end{proposition}

{\bf Proof:}
From (\cite{Hal},9), ou (\cite{FHT},14), the morphism $f$ has a Sullivan minimal model, ie, there exists a commutative diagram of CDGA algebras:
$$\begin{array}{ccc}
(A,d_A)&\stackrel{f}{\longrightarrow }&(B,d_B)\\
{\scriptstyle m_A} \uparrow && \uparrow {\scriptstyle m_B}\\
(\bigwedge U,d)&\stackrel {\scriptstyle i}{\longrightarrow }& (\bigwedge U\otimes \bigwedge V,d')
\end{array}$$
where $(\bigwedge U,d)$ is a Sullivan minimal algebra,  $(\bigwedge U,d) \longrightarrow  (\bigwedge U\otimes \bigwedge V,d')$ is a minimal relative Sullivan algebra and the vertical maps are quasi-isomorphisms. The existence of $g$ satisfying $g\circ f=Id$ implies that there exists a retraction $q: (\bigwedge U\otimes \bigwedge V,d')\longrightarrow  (\bigwedge U,d)$ satisfying $q\circ i$ homotopic to the identity map. A classical argument  implies that $d'$ is minimal. Now we use the fact that $(B,d_B)$ is formal, so there exists a CDGA map $\rho : (\bigwedge U\otimes \bigwedge V,d')\longrightarrow (H^*(B),0)$
 such that $\rho ^*=m_B^*$. Consider
 $$\theta =g^*\circ \rho \circ i: (\bigwedge U,d)\longrightarrow (H^*(A),0)$$ 
Then we have: 
$\theta ^*=g^*\circ \rho ^*\circ i*=g^*\circ m_B^*\circ i^*=g^*\circ f^*\circ m_A^*=m_A^*$. This proves that $(A,d_A)$ is formal.

\bigskip

\section{Proof of the Main Theorem}
It relies on the results of \cite{Si} and \cite{Vi} and does not use the computations of \cite{Hae} or \cite{BS}. 

It is an immediate consequence of theorem 3.1 and theorem 3.5.

\begin{thm}
Let $X$ and $Y$ be  spaces satisfying the hypothesis of the Main Theorem. If $\mathcal F(X,Y)$ is formal then there exists an integer $2d+1\geq 1$  such that  $\mathcal F(S^{2d+1},Y)$ is formal.
\end{thm}

To prove Theorem 3.1, we will use the Lie model for the space $\mathcal F(X,Y)$ explained in  Proposition 3.2.

Since $X$ is a nilpotent space with finite dimensional cohomology, it has a finite dimensional model $(A,d_A)$ in the category of commutative differential graded algebras (Proposition 2.1).

\begin{proposition}(\cite{Si},section 6), (\cite{Vi} theoreme 1).
Let $X$ be a nilpotent space
 such that there exists $N\geq 1$ with $H^N(X,\mathbb Q) \ne 0$ and $H^n(X,\mathbb Q)=0$ if $n>N$. Let $Y$ be a m-connected 
space with $m\geq N+1$. If $(A,d_A)$ is a finite dimensional model of $X$ satisfying $A^n=0$ if $n>N$. If $(L,d_L)$ is a Lie model of $Y$, then $(A\otimes L,D)$ is a Lie model for the space $\mathcal F(X, Y)$ where the structure of differential graded Lie algebra on $A\otimes L$ is the following: 
\begin{enumerate}
\item $\vert a\otimes l\vert=-\vert a\vert +\vert l \vert$ if  $a\in A$,\,\, $l\in L$
\item $ [a\otimes l,a'\otimes l']=(-1)^{\vert a'\vert \cdot \vert l\vert} aa'\otimes [l,l']$
\item $D(a\otimes l)=d_Aa\otimes l+(-1)^{\vert a\vert }a\otimes d_L(l)$
\end{enumerate}
Furthermore the projection:  $(A,d_A) \mapsto A^0=\mathbb Q$ extends to a morphism of differential Lie algebras: $(A\otimes L,D)\rightarrow (L,d_L)$ which is a model of the fibration $p: \mathcal F(X,Y)\rightarrow Y,$ defined by $p(f) =f(x_0)$ if $x_0$ is a fixed point in $X$ and the inclusion $\mathbb Q \hookrightarrow A$ extends to a morphism of differential Lie algebras which is a model of the canonical section of the fibration $p$.

\end{proposition}

This Lie model  permits us to prove two lemmas (lemma 3.3 and lemma 3.4), the first one is used to prove theorem 3.1 and the second one  is used to prove theorem 3.5.

\begin{lem} Let $X$ be a space satisfying the hypothesis of the main theorem. Let $(A,d_A)$ be a finite dimensional model of $X$. Then there exists an exterior algebra $\bigwedge t$ with $\vert t\vert =2d+1\geq 1$ and morphisms $i: (\bigwedge t,0)\rightarrow (A,d_A)$, $q:(A,d_A)\rightarrow (\bigwedge t,0)$ in the CDGA category satisfying: $q\circ i=Id$.
\end{lem}

{\bf Proof of lemma 3.3.} 
If $X$ has a minimal Sullivan model $(\bigwedge V,d)$ with ${\rm Ker} d\cap V_{odd}\ne0$.  Let $t$ be an odd generator of $V$ with $dt=0$ and $\vert t\vert=2d+1$. We denote by $i_0$ the inclusion $(\bigwedge t,0)\rightarrow (\bigwedge V,d)$, the linear projection $V\mapsto \mathbb Q t$ extends to a morphism of differential algebras $q_0: (\bigwedge V,d)\rightarrow (\bigwedge t,0)$ since $d$ is minimal. Let $m: (\bigwedge V,d)\rightarrow (A,d_A)$  be the finite dimensional model of $X$ , (Proposition 2.1). From the construction of $(A,d_A)$, it is easy to check that $q_0$ factors through $(A,d_A$),ie,  there exists $q: (A,d_A)\rightarrow (\bigwedge t,0)$ such that $q\circ m=q_0$. Put $i=m\circ i_0$ then we have: $q\circ i= Id$.

\medskip
 {\bf Proof of theorem 3.1}. Using  Proposition 3.2 and Lemma 3.3, we  define differential Lie morphisms $I=i\otimes _{Id}:(\bigwedge t,0)\otimes (L,d_L) \longrightarrow (A,d_A)\otimes (L,d_L)$ and $Q=q\otimes Id:   (A,d_A)\otimes (L,d_L)\longrightarrow (\bigwedge t,0)\otimes (L,d_L)$ satisfying $Q\circ I=Id$. Here $(A,d_A)$ is a finite dimensional model of $X$ in the category CDGA, $(L,d_l)$ is a Lie model of $Y$ and $(\bigwedge t,0)$ is a finite model of $S^{2d+1}$. Let $\mathcal C^*$ be the functor defined in section 2 from the category of  differential Lie algebras to the category of commutative differential algebras. Denote $g=\mathcal C^*(I): \mathcal C^*(A\otimes L) \longrightarrow \mathcal C(\bigwedge t \otimes L)$ and $f=\mathcal C^*(Q): \mathcal C(\bigwedge t \otimes L)\longrightarrow \mathcal C^*(A\otimes L) $, $f$ and $g$ are two morphisms in the category CDGA, we use Proposition 2.2 to conclude.

\begin{lem}
Let $X$ be a nilpotent space
 such that there exists $N\geq 1$ with $H^N(X,\mathbb Q) \ne 0$ and $H^n(X,\mathbb Q)=0$ if $n>N$. Let $Y$ be a m-connected 
space with $m\geq N+1$. If $\mathcal F(X,Y)$ is formal, then $Y$ is formal.
\end{lem}

{\bf Proof of lemma 3.4}. From Proposition 3.2, the projection:  $(A,d_A) \mapsto A^0=\mathbb Q$ extends to a morphism of differential Lie algebras: $q:(A\otimes L,D)\rightarrow (L,d_L)$ which is a model of the fibration $p: \mathcal F(X,Y)\rightarrow Y,$ the inclusion $\mathbb Q \hookrightarrow A$ extends to a morphism of differential Lie algebras $i:(L,d_L) \hookrightarrow (A,d_A)\otimes (L,d_L)$   which is a model of the canonical section of the fibration $p$, so that $q\circ i=Id$. Let $\mathcal C^*$ be the functor defined in section 2 from the category of  differential Lie algebras to the category of commutative differential algebras. Denote $g=\mathcal C^*(i)$ and $f=\mathcal C^*(q)$. As above, we use Proposition 2.2 to conclude.

This lemma is proved also in \cite{DV} and \cite{Y}, using other arguments.

\medskip

Now we come back to the category CDGA and we will prove the following theorem 
 using similar arguments to those developed in \cite{DV}.

\begin{thm}
Let $Y$ be a m-connected space  such that  $H^*(Y,\mathbb Q)$ is finitely generated as algebra. 
We assume that there exists  $p\geq 1$, $p$ odd, with $m\geq p+1$ such that $\mathcal F(S^{p},Y)$  is formal. Then
$Y$ has the rational homotopy type of a finite product of Eilenberg Mac Lane spaces.
\end{thm}

{\bf Proof of Theorem 3.5}.
We have proved in Lemma 3.4 that $Y$ is formal.  Put $\bigwedge t=H^*(S^p)$.  Let $(L,d_L)$ be a Lie model of $Y$, then $\mathcal C^*(L)=(\bigwedge V,d)$ is a Sullivan algebra. 
From proposition 3.2, a Sullivan model of $\mathcal F(S^p,Y)$ is $\mathcal C^*(L\oplus \bar L,D)$, where $\bar L_n=\mathbb Qt\otimes L_{n+p} \simeq L_{n+p}$, $D_{\vert L}=d_L$,  $D\bar x=-\overline{d_L(x)}$, $(-1)^{\vert a\vert}[a,\bar b]=\overline{[a,b]}$, $[\bar a,\bar b]=0$. So we have $\mathcal C^*(L\oplus \bar L,D)=(\bigwedge(V\oplus SV),d)$ with $SV^n=V^{n+p}$.   The inclusion $\mathcal C^*(L)\hookrightarrow \mathcal C^*(L\oplus \bar L,D)$ is a relative Sullivan model of the fibration  $\mathcal F(S^p,Y)\rightarrow Y.$
  We extend the identity map  $S: V\rightarrow SV$ to a derivation of graded algebras
 of degree $-p$: $\bigwedge V\rightarrow \bigwedge V\otimes \bigwedge SV$. The differential $d$ on  $\mathcal C^*(L\oplus \bar L,D)=\bigwedge(V\oplus SV)$ is defined by the condition $d(Sv)=-S(dv)$ for $v\in V$. \\
Since $Y$ is formal, we will work with its bigraded minimal model $(\bigwedge Z,d)$ in the sense of Halperin-Stasheff, \cite{HS}, $Z=\oplus_n Z^n$,  $Z^n=\oplus _{k\geq 0}Z^n_k$ and $d(Z^n_k)\subset (\bigwedge Z )^{n+1}_{k-1}$. Since $H^*(Y,\mathbb Q)$ is finitely generated as algebra, we have $\dim Z_0<\infty $ where $Z_0=\oplus_n Z_0^n$. 
 Let $\bar Z^n=Z^{n+p}$ and $\bar Z=\oplus Z^n$, the identity map $S$ defined by $S(z)=\bar z$ is extended to an algebra derivation of degree $-p$  on $(\bigwedge (Z\oplus \bar Z),d)$ by putting $S(\bar z)=0.$ It is clear that $(\bigwedge (Z\oplus \bar Z),d)$ is a minimal Sullivan model of $\mathcal F(S^p,Y)$ where $d(\bar z)=-S(dz)$.

A  generalization of lemme 3 in \cite{DV} can be formulated as follows.

\begin{proposition}
Let $X=S^p$, $p$ odd, and $Y$ be spaces satisfying the hypothesis of the main theorem. If $\mathcal F(S^p,Y)$ is formal and $(\bigwedge Z,d)$ is the bigraded minimal model  of $Y$, then the minimal Sullivan algebra $(\bigwedge (Z\oplus \bar Z),d)$ bigraded by $(\bar Z)_n=\overline{(Z_n)}$ is the bigraded model of $\mathcal F(S^p,Y)$ in the sense of Halperin-Stasheff.
\end{proposition}

The proof of this proposition is postponed to the end of this section.
\medskip

Recall a key lemma (lemme 1) in \cite{DV}.

\begin{lem}
 Let $(\bigwedge(W_0\oplus W_+),d)$ be the bigraded minimal model of a formal space such that $\dim W_0<\infty $. Then for any nonzero element in $W_+^{even}$, there exist an element $w'\in W_+^{odd}$, an integer $n\geq 2$, and a decomposable element $\Omega$ without nonzero component in $w^n$ such that $dw'=w^n+\Omega $.
\end{lem}  

Now we finish the proof of Theorem 3.5.

Consider the bigraded model of  $\mathcal F(S^p,Y)$ defined by proposition 3.6. Suppose that $\bar Z^{even}\ne 0$ and consider a nonzero element $  \bar z \in \bar Z^{even}$. Since $d(\bar Z)\subset \bar Z\cdot \bigwedge Z$, there does not exist element $w'\in Z\oplus \bar Z$ such that $dw'=\bar z^n+\Omega $ for some $n\geq 2$. Lemma 3.7 implies that $\bar Z_+^{even}=0$. Since $p$ is odd, and $\bar Z^n=Z^{n+p}$, it follows  that $Z_+^{odd}=0$. If we apply lemma 3.7 to $(\bigwedge Z,d)$ we get $Z_+^{even}=0$. Finally we have $Z=Z_0$ and $d=0.$\\
This achieves the proof of  Theorem 3.5 and also the proof of the Main Theorem.\\
 We note that if $p$ was even we could not conclude.

\medskip

{\bf Proof of Proposition 3.6}. We have to prove that $H_n(\bigwedge (Z\oplus \bar Z),d)=0$ for all $n\geq 1$. Suppose that this affirmation is not true. Let $q$ be the lowest non-zero integer such that $H_q(\bigwedge (Z\oplus \bar Z),d)\ne 0$. Let $l$ be the lowest non-zero integer such that $H_q^l(\bigwedge (Z\oplus \bar Z),d)\ne 0$. Let $[\alpha ]$ be a non-zero element in $H_q^l(\bigwedge (Z\oplus \bar Z),d)$. It is easy to check that $[\alpha ]$ does not belong to $H^+\cdot H^+$ where $H^+={\displaystyle \sum _{n>0}H^n(\bigwedge (Z\oplus \bar Z),d)}$. Since $(\bigwedge (Z\oplus \bar Z),d)$ is the minimal model of a formal space, it is classical to prove that there exist $z\in Z_q$ and $\gamma \in (\bigwedge Z \otimes  \bar Z)_q$ such that $\alpha =\bar z + \gamma $, \cite {DV}. It remains to prove that such a cocycle cannot occur in $(\bigwedge Z \otimes  \bar Z)_q$ with $q\geq 1$. \\
In \cite{DV}, Proposition 3.6 is proved when $p=1$ using technical calculations. Here we give a direct proof for any odd integer $p$.\\
The map  $\tilde S: (\bigwedge (Z\oplus \bar Z),d)\rightarrow (\bigwedge (Z\oplus \bar Z),d)$ defined by $S(a)=(-1)^{\vert a \vert} S(a)$ for $a\in \bigwedge (Z\oplus \bar Z)$ is a morphism of cochain complexes of degree $-p$. We have a short exact sequence of complexes:
$$0 \rightarrow  (Ker S,d) \xrightarrow{j} (\bigwedge (Z\oplus \bar Z),d)\xrightarrow{\tilde S} (Im S,d)_{-p} \rightarrow 0$$
where $j$ is the inclusion.\\
Furthermore, we have: $(Ker S)^n= (Im S)^n$ for all $n\geq 1$. Since $(\bigwedge (Z\oplus \bar Z),d)$ is formal, a reformulation of theorem A of \cite{Vi2} in this context implies that the cohomology long exact sequence associated to the exact sequence above splits into short exact sequences:
$$0 \rightarrow  H^{n+p}(Ker S,d) \rightarrow H^{n+p}(\bigwedge (Z\oplus \bar Z),d)\rightarrow H^n(Im S,d) \rightarrow 0$$
for all $n\geq 1$.\\
Now we work with the cocycle $\alpha =\bar z+\gamma $, we have $0=d\alpha =dS(z) +d\gamma $. Since $S^2=0$, we get $d(S(\gamma )=-Sd\gamma =0$. So $S\gamma $ is a cocycle in $Im S$. Since $\tilde S^*$ is surjective in cohomology, there exists a cocycle $\beta \in (\bigwedge (Z\oplus \bar Z))^{l+p}$ so that $[S\beta ]=[S\gamma ]$ in $H^*(Im S,d)$. So there exists $\mu$ such that $S\beta =S\gamma +dS\mu $, that is $S(\gamma -\beta -d\mu )=0.$ Since $Ker S= Im S$, there exists $\varphi \in (\bigwedge Z)_q$ such that $\gamma -\beta -d\mu =S\varphi $ and $\varphi $ is decomposable.\\
Recall that $\alpha =\bar z +\gamma $ is a cocycle in $(\bigwedge Z \otimes  \bar Z)_q$ with $q\geq 1$. We have $\alpha = S(z+\varphi )+ \beta +d\mu $. Put $z'=z+\varphi $ and $\beta '=\beta +d\mu $, we have $\alpha = \bar z' +\beta '$, so $d\bar z'=0$. Since $z'\in (\bigwedge^+ Z)_q$, we get $dz'=0$. This is a contradiction with the fact that $(\bigwedge Z,d)$ is the bigraded minimal model in the sense of \cite{HS}. 
\bigskip

\section{A counterexample to the non-formality of $\mathcal F(S^p,Y)$ when $p$ is even.}
Let $Y= K(\mathbb Q,4)\vee K(\mathbb Q,4) $, it is  a $3$-connected space whose minimal Sullivan model is $(\bigwedge(x_1, x_2, y),d)$ with $dx_1=0,\, \, dx_2=0,\,\, dy=x_1x_2$, $\vert x_1\vert=\vert x_2\vert=4$ and $\vert y\vert=7$. We have $H^*(Y,\mathbb Q)=\mathbb Q[x_1,x_2]/(x_1x_2)$ and $Y$ is formal.
The propositions proved in section 3 show that a minimal model of $\mathcal F(S^2,Y)$ is the following:
$$(\bigwedge,d)=(\bigwedge (x_1, x_2, y, \bar x_1, \bar x_2, \bar y),d)$$
with $\vert \bar x_1\vert = \vert \bar x_2\vert = 2$ and $\vert \bar y\vert= 5$. We have $d\bar x_1=d\bar x_2=0$ and $d\bar y=\bar x_1x_2+x_1\bar x_2$. It is easy to check that the polynomials $(dy, d\bar y)$ form a  regular sequence in $\mathbb Q[x_1, x_2]$ so $(\bigwedge,d)$ is a Koszul complex, hence it is formal.

 \vspace{1cm}
\hspace{-1cm}\begin{minipage}{19cm}
 \small
vigue@math.univ-paris13.fr\\ D\'epartement de math\'ematiques  
\\
 Institut
Galil\'ee, \\
Universit\'e de
Paris-Nord\\
 93430
Villetaneuse, France
\end{minipage}

\end{document}